\documentclass{article}

\usepackage{amsmath,amsfonts,amsthm,amssymb,graphicx,tikz,tikz-cd}
\usepackage[all,2cell,ps]{xy}

\theoremstyle{plain}
\newtheorem{thm}{Theorem}[section]
\newtheorem{lem}[thm]{Lemma}
\newtheorem{prop}[thm]{Proposition}

\theoremstyle{definition}

\newcommand{\Z}{\mathbb Z}

\newcommand{\Q}{\mathbb Q}
\newcommand{\B}{\mathbb B}
\newcommand{\C}{\mathbb C}

\newcommand{\sig}{\sigma}
\newcommand{\al}{\alpha}

\newcommand{\Gam}{\Gamma}
\newcommand{\del}{\delta}

\DeclareMathOperator{\PU}{PU}

\newcommand{\ssm}{\smallsetminus}
\newenvironment{pf}{\begin{proof}}{\end{proof}}

\usetikzlibrary{arrows}

\title{Multiple realizations of varieties as ball quotient compactifications}
\author{Luca F.\ Di Cerbo \\ \small{Notre Dame University} \\ \small{\textsf{ldicerbo@nd.edu}} \and Matthew Stover\footnote{This material is based upon work supported by the National Science Foundation 
under Grant Number NSF DMS-1361000. The second author acknowledges support from U.S.\ National Science Foundation grants DMS 1107452, 1107263, 1107367 ``RNMS: GEometric structures And Representation varieties'' (the GEAR Network).} \\ \small{Temple University}\\ \small{\textsf{mstover@temple.edu}}}
\date{\today}

\begin{document}

\maketitle

%%%%%%%%%%%%%%%%%%%%
\begin{abstract}
We study the number of distinct ways in which a smooth projective surface $X$ can be realized as a smooth toroidal compactification of a ball quotient. It follows from work of Hirzebruch that there are infinitely many distinct ball quotients with birational smooth toroidal compactifications. We take this to its natural extreme by constructing arbitrarily large families of distinct ball quotients with biholomorphic smooth toroidal compactifications.
\end{abstract}
%%%%%%%%%%%%%%%%%%%%

%%%%%%%%%%%%%%%%%%%%
\section{Introduction}\label{sec:Intro}
%%%%%%%%%%%%%%%%%%%%

Let $\B^2$ be the unit ball in $\C^2$ with its Bergman metric and $\Gam \subset \PU(2,1)$ a nonuniform lattice. Then $Y = \B^2 / \Gam$ is a complex orbifold with a finite number of cusps, and it admits a number of compactifications by a normal projective variety, which may or may not be smooth. When $\Gam$ is neat (see \S \ref{ssec:Ball}), then $Y$ has a particularly nice \emph{toroidal compactification}, which is a smooth projective surface \cite{AMRT, Mok}. It is an important and open question to decide which projective surfaces are compactifications of ball quotients, and more specifically which smooth projective surfaces are smooth toroidal compactifications.

In this paper, we study the number of ways in which a fixed projective surface can arise as a compactification of a ball quotient. A by-product of Hirzebruch's construction of smooth projective surfaces with $c_1^2 / c_2$ arbitrarily close to $3$ \cite{Hirzebruch} is an infinite family of ball quotients with \emph{birational} smooth toroidal compactifications (cf.\ \cite{Holzapfel}, and see \cite{Kasparian} for more birational examples). The purpose of this paper is to exhibit arbitrarily large families of distinct ball quotients with \emph{biholomorphic} smooth toroidal compactifications.

%%%%%%%%%%%%%%%%%%%%
\begin{thm}\label{thm:MainThm}
For any natural number $n$, there exists a smooth projective surface $X = X(n)$ and pairwise nonisomorphic quotients $Y_1, \dots, Y_n$ of $\B^2$ by neat lattices $\Gam_1, \dots, \Gam_n$ in $\PU(2, 1)$ such that $X$ is a smooth toroidal compactification of $Y_i$ for each $1 \le i \le n$.
\end{thm}
%%%%%%%%%%%%%%%%%%%%

We note that finiteness is necessary, i.e., a fixed smooth projective surface can only arise as a smooth toroidal compactification in finitely many ways. Indeed, if $Y = \B^2 / \Gam$ is a ball quotient with smooth toroidal compactification $X$, then $Y$ and $X$ have the same topological Euler number. Since the volume of $Y$ in the metric descending from $\B^2$ is proportional to its Euler number \cite{Mumford}, all ball quotients with smooth toroidal compactification $X$ have the same volume. However, H.\ C.\ Wang showed that there are only finitely many isomorphism classes of lattices in $\PU(2,1)$ of bounded covolume \cite{Wang}, and the claim follows.

Another interpretation of Theorem \ref{thm:MainThm} is the following. Recall that if $X$ is a smooth toroidal compactification of $Y_i$ then there exists a divisor $D_i$ on $X$, where $D_i$ is a disjoint union of elliptic curves with negative self-intersection, such that $Y_i = X \ssm D_i$. Equivalently, if a pair $(X, D_i)$ saturates (i.e., realizes equality in) the logarithmic Bogomolov--Miyaoka--Yau inequality
\[
\overline{c}_1^2(X, D_i) \le 3 \overline{c}_2(X, D_i),
\]
then $Y_i = X \ssm D_i$ is a ball quotient (see \S \ref{ssec:Ball}). Recall that $\overline{c}_1^2$ denotes the self-intersection of the log canonical divisor $K_X + D_i$ and $\overline{c}_2$ is the topological Euler number of $X \ssm D_i$. In Theorem \ref{thm:MainThm} we then produce, for any given positive integer $n$, a smooth surface $X=X(n)$ and divisors $D_1, \dots, D_n$ on $X$ such that each $X \ssm D_i$ is a distinct quotient of the ball by a neat lattice.

The paper is organized as follows. Section \ref{sec:Prelim} collects preliminary facts about the ball and its Bergman metric, smooth toroidal compactifications, and a ball quotient of Euler number $1$ constructed by Hirzebruch in \cite{Hirzebruch}. This particular example plays a fundamental role in our construction. Finally in \S \ref{sec:Proof}, we  prove Theorem \ref{thm:MainThm}.\\

\noindent\textbf{Acknowledgements}. We thank the organizers and participants of the conference ``Algebraic \& Hyperbolic Geometry - New Connections'' for the stimulating environment that set this collaboration in motion.

%%%%%%%%%%%%%%%%%%%%
\section{Preliminaries}\label{sec:Prelim}
%%%%%%%%%%%%%%%%%%%%

%%%%%%%%%%%%%%%%%%%%
\subsection{Ball quotients and their compactifications}\label{ssec:Ball}
%%%%%%%%%%%%%%%%%%%%

Let $\B^2$ be the unit ball in $\C^2$ with its Bergman metric. See \cite{Goldman} for more on its geometry.  The group of biholomorphic isometries of $\B^2$ is isomorphic to $\PU(2, 1)$, and a discrete subgroup $\Gam \subset \PU(2, 1)$ is a \emph{lattice} if $Y = \B^2 / \Gam$ has finite volume in the metric descending from $\B^2$. Then $Y$ is a manifold if and only if $\Gam$ is torsion free. When $Y$ is not compact, it is well-known that it is a quasiprojective variety and has a finite number of topological ends, or \emph{cusps}.

Suppose further that $\Gam$ is \emph{neat}, i.e., that the subgroup of $\C$ generated by the eigenvalues of $\Gam$ is torsion free. In particular, $\Gam$ is torsion free. This implies that $Y$ admits a smooth toroidal compactification $Y^*$ by adding a certain elliptic curve at each cusp. (Actually, one only needs that the subgroup of $\C$ generated by the eigenvalues of parabolic elements of $\Gam$ is torsion free, which means that the parabolic elements have no rotational part.) See \cite{AMRT} and \cite{Mok} for more details.

Let $Y$ be the quotient of $\B^2$ by a neat lattice and $Y^*$ its smooth toroidal compactification. Then $Y^* \ssm Y$ consists of a finite union of disjoint elliptic curves $T_i$, each having negative self-intersection. Let $D = \sum T_i$. Hirzebruch--Mumford proportionality \cite{Mumford} implies that
\begin{equation}\label{eq:HMp}
\overline{c}_1^2(Y^*, D) = 3 \overline{c}_2(Y^*, D).
\end{equation}
Moreover, the logarithmic version of Yau's solution to the Calabi conjecture implies the converse \cite{Tian--Yau}. More precisely, if $X$ is a smooth projective variety and $D$ a normal crossings divisor on $X$ satisfying \eqref{eq:HMp}, then $X$ is the smooth toroidal compactification of
\[
Y = X \ssm D = \B^2 / \Gam
\]
for some torsion free lattice $\Gam \subset \PU(2, 1)$.

%%%%%%%%%%%%%%%%%%%%
\subsection{Hirzebruch's example}\label{ssec:Hirz}
%%%%%%%%%%%%%%%%%%%%

We now describe an example from \cite{Hirzebruch}, which is critical in our construction. Let $\zeta = e^{\pi i / 3}$ and $\rho = \zeta^2$. Then $\Z[\rho]$ is a lattice in $\C$, and $E = \C / \Z[\rho]$ is the elliptic curve of $j$-invariant $0$. Set $A = E \times E$, let $[z, w]$ be coordinates on $A$, and consider the curves:
\begin{align*}
T_0 &= \{w = 0\} \\
T_\infty &= \{z = 0\} \\
T_1 &= \{w = z\} \\
T_\zeta &= \{w = \zeta z\}
\end{align*}

It is easy to check that each $T_\al$ is isomorphic to $E$, has self-intersection $0$, and that $T_\al \cap T_\beta = \{[0,0]\}$ for each $\al \neq \beta$. Let $A^*$ be the blowup of $A$ at $[0,0]$ and $D$ the proper transform of $\sum T_\al$ in $A^*$. It is easy to check that $\overline{c}_1^2(A^*, D) = 3 \overline{c}_2(A^*, D)$. In particular, $Y = A^* \ssm D$ is the quotient of the ball by a torsion free lattice and $A^*$ is the smooth toroidal compactification of $Y$.

%%%%%%%%%%%%%%%%%%%%
\section{Proof of Theorem \ref{thm:MainThm}}\label{sec:Proof}
%%%%%%%%%%%%%%%%%%%%

Let $E$ and $A = E \times E$ be as in \S\ref{ssec:Hirz}. The first step in our construction is to exploit certain self-isogenies of $A$. Translation by a generator for the $\Q$-rational $3$-torsion subgroup $E(\Q)[3] \cong \Z / 3$ of $E$ determines a degree $3$ self-isogeny $r : E \to E$. We recast this in the language of covering space theory. Let $v_1 = 1$ and $v_2 = \rho$ be the usual generators for $\Z[\rho] = \pi_1(E)$. We can realize $r^n$ via the $3^n$-fold covering of $E \to E$ associated with the kernel of the homomorphism
\[
R_j : \pi_1(E) \to \Z / 3^n \quad\quad R_j(v_1) = R_j(v_2) = \del,
\]
where $\del$ is a chosen generator for $\Z / 3^n$. Notice that this kernel is the subgroup $(1 - \rho)^n \Z[\rho]$, which is the unique ideal of norm $3^n$ in the ring $\Z[\rho]$.

Now, let $A = E \times E$ and
\begin{align*}
&v_1 = \begin{pmatrix} 1 \\ 0 \end{pmatrix} &v_3 = \begin{pmatrix} 0 \\ 1 \end{pmatrix} \\
&v_2 = \begin{pmatrix} \rho \\ 0 \end{pmatrix} &v_4 = \begin{pmatrix} 0 \\ \rho \end{pmatrix}
\end{align*}
be generators for $\pi_1(A)$. Fix $n$ and let $\del$ be a generator for $\Z / 3^n$. For any $0 \le j < 3^n$, we obtain a homomorphism $\sig_{n, j} : \pi_1(A) \to \Z / 3^n$ by:
\begin{align*}
\sig_{n, j}(v_1) = \sig_{n, j}(v_2) &= \del \\
\sig_{n, j}(v_3) = \sig_{n, j}(v_4) &= \del^j
\end{align*}
This induces a self-isogeny $s_{n, j} : A \to A$ that restricts to an isogeny of degree $3^n$ on the first factor and of degree $|\del^j| = 3^n / \gcd(j, 3^n)$ on the second. Note that the kernel $K_{n, j}$ of $\sig_{n, j}$ is generated by $3^n v_1$, $v_1 - v_2$, $v_3 - v_4$, and $j v_1 - v_3$.

Let $A^*$ be the blowup of $A$ at the origin $[0,0]$. Then $\sig_{n, j}$ induces an \'etale covering $s_{n, j}^* : A^*_{n, j} \to A^*$, where $A_{n, j}^*$ is the blowup of $A$ at $s_{n, j}^{-1}([0,0])$. Let $Y \subset A^*$ be Hirzebruch's ball quotient. We now need to know that $s_{n, j}^*$ induces a connected covering of $Y$.

%%%%%%%%%%%%%%%%%%%%
\begin{lem}\label{lem:Covers}
Let $\pi : B \to A$ be a finite \'etale cover with group $G$. Then there exists a connected ball quotient $Y^\prime$ and a finite regular covering $Y^\prime \to Y$ induced by $\pi$ with Galois group $G$.
\end{lem}
%%%%%%%%%%%%%%%%%%%%

%%%%%%%%%%%%%%%%%%%%
\begin{pf}
Since $\pi_1(A^*) \cong \pi_1(A)$, we also obtain an \'etale covering $\pi^* : B^* \to A^*$, where $B^*$ is the blowup of $B$ above each of the $\deg(\pi)$ points in $\pi^{-1}([0,0])$. Let $Y^\prime$ be the inverse image in $B^*$ of $Y$ under $\pi^*$. Then $Y^\prime \to Y$ is an \'etale cover. We must show that $Y^\prime$ is connected and that the Galois group of this covering is $G$. However, it is well known that the map
\[
\pi_1(Y) \to \pi_1(A^*) \cong \pi_1(A)
\]
induced by the inclusion is onto \cite[Prop.\ 2.10]{Kollar}. In particular, the induced homomorphism $\pi_1(Y) \to G$ associated with the covering $B^* \to A^*$ is onto. The lemma follows from elementary covering space theory.
\end{pf}
%%%%%%%%%%%%%%%%%%%%

It follows immediately from Lemma \ref{lem:Covers} that the self-isogenies $s_{n, j}$ induce \'etale covers $Y_{n, j} \to Y$ with degree $3^n$. Moreover, it is clear from the construction that $Y_{n, j}$ is a quotient of the ball by a neat lattice and that the smooth toroidal compactification is the blowup $A$ at $3^n$ distinct points. We now count the number of cusps of $Y_{n, j}$.

%%%%%%%%%%%%%%%%%%%%
\begin{prop}\label{prop:CuspCount}
For and $n \ge 0$ and $0 \le j < 3^n$, the ball quotient $Y_{n, j}$ has
\[
\begin{cases}
6 & j \equiv 1 \pmod{3} \\
3 + \gcd(j + 1, 3^n) & j \equiv 2 \pmod{3} \\
3 + \gcd(j, 3^n) & j \equiv 0 \pmod{3}
\end{cases}
\]
cusps.
\end{prop}
%%%%%%%%%%%%%%%%%%%%

%%%%%%%%%%%%%%%%%%%%
\begin{pf}
We need to count the number of lifts of each $T_\al$, $\al \in \{0, \infty, 1, \zeta\}$, under $s_{n, j}$. Equivalently, we need to calculate the index of $\pi_1(T_\al)$ in $\Z / 3^n$ under $\sig_{n, j}$. As subgroups of $\pi_1(A)$, we have:
\begin{align*}
\pi_1(T_0) &= \langle v_1, v_2 \rangle&  \pi_1(T_\infty) &= \langle v_3, v_4 \rangle \\
\pi_1(T_1) &= \langle v_1 + v_3, v_2 + v_4\rangle&  \pi_1(T_\zeta) &= \langle v_1 + v_3 + v_4, v_2 - v_3 \rangle \\
\end{align*}
Therefore:
\begin{align*}
\sig_{n, j} \big( \pi_1(T_0) \big) &= \Z / 3^n  &&\sig_{n, j} \big( \pi_1(T_\infty) \big) = \langle \del^j \rangle \\
\sig_{n, j} \big( \pi_1(T_1) \big) &= \langle \del^{j + 1}\rangle &\Z / 3^{n - 1} \subseteq \ &\sig_{n, j} \big( \pi_1(T_\zeta) \big) \ \subseteq \Z / 3^n
\end{align*}
The only nontrivial point is the last one. The image certainly contains $\del^{1 + 2 j}$ and $\del^{1 - j}$, the images of the two generators. If $j$ is congruent to $0$ or $2$ modulo $3$, then either of these generates all of $\Z / 3^n$. However, it always contains
\[
\del = \del^{1 + 2 j} \del^{2(1 - j)} = \del^3,
\]
which generates $\Z / 3^{n - 1}$. When $j$ is congruent to $1$ modulo $3$, it follows that the image is $\Z / 3^{n - 1}$.

There is always exactly one cusp above $T_0$. There are exactly $3$ cusps above $T_\zeta$ when $j \equiv 1 \pmod{3}$, and there is exactly $1$ cusp otherwise. For $T_\infty$ and $T_1$, we must calculate the index of $\langle \del^j \rangle$ and $\langle \del^{j + 1} \rangle$ in $\Z / 3^n$, respectively. Note that at least one of $j$ and $j + 1$ is not divisible by $3$, hence one of $\del^j$ or $\del^{j + 1}$ generates all of $\Z / 3^n$. Since $\del^k$ has order $3^n / \gcd(3^n, k)$ for any $0 \le k < 3^n$, the proposition follows.
\end{pf}
%%%%%%%%%%%%%%%%%%%%

Now we prove that the spaces $Y_{n, j}$ all have isomorphic smooth toroidal compactifications.

%%%%%%%%%%%%%%%%%%%%
\begin{prop}\label{prop:IsomOfSTC}
For any fixed $n$ and $0 \le j, k < 3^n$, the spaces $A_{n, j}^*$ and $A_{n, k}^*$ are isomorphic.
\end{prop}
%%%%%%%%%%%%%%%%%%%%

%%%%%%%%%%%%%%%%%%%%
\begin{pf}
Let $r$ be an integer and consider the linear map $f : \C^2 \to \C^2$ with matrix
\[
\begin{pmatrix} 1 & r \\ 0 & 1 \end{pmatrix}.
\]
Then, in terms of the $\Z$-module generators $\{v_i\}$ of $\pi_1(A)$, we have:
\begin{align*}
f(v_1) &= v_1 \\
f(v_2) &= v_2 \\
f(v_3) &= r v_1 + v_3 \\
f(v_4) &= r v_2 + v_4
\end{align*}
Now, we consider what $f$ does to the kernel $K_{n, j}$ of $\sig_{n, j}$. Since
\begin{align*}
f(3^n v_1) &= 3^n v_1 \\
f(v_1 - v_2) &= v_1 - v_2 \\
f(v_3 - v_4) &= r (v_1 - v_2) + (v_3 - v_4) \\
f(j v_1 - v_3) &= (j - r) v_1 - v_3
\end{align*}
we see that the image of $K_{n, j}$ also contains
\[
v_3 - v_4 = f(v_3 - v_4) - r f(v_1 - v_2).
\]
If $0 \le k < 3^n$ is the representative for $j - r$ modulo $3^n$, it follows that $f$ contains generators for the kernel $K_{n, k}$ of $\sig_{n, k}$. Since the matrix of $f$ has determinant one, it follows that $f(K_{n, j}) = K_{n, k}$.

However, $f$ induces an isomorphism of  $\Z[\rho]\times \Z[\rho]$ onto itself. Thus, $f$ induces a commutative diagram of isomorphisms and coverings:
\[
\begin{tikzcd}
A \arrow{r} & A \\
A_{n, j} \arrow{u} \arrow{r} & A_{n, k} \arrow{u}
\end{tikzcd}
\]
taking $[0,0]$ to itself in $A$. In particular, the map $A_{n, j} \to A_{n, k}$ takes $s_{n, j}^{-1}([0,0])$ to $s_{n, k}^{-1}([0,0])$. This proves that $f$ induces an isomorphism $A_{n, j}^* \to A_{n, k}^*$.
\end{pf}
%%%%%%%%%%%%%%%%%%%%

Now, we have the tools necessary to prove Theorem \ref{thm:MainThm}.

%%%%%%%%%%%%%%%%%%%%
\begin{pf}[Proof of Theorem \ref{thm:MainThm}]
Fix $n$. Then we can find $0 \le j_1, \dots, j_n < 3^n$ for which Proposition \ref{prop:CuspCount} implies that the ball quotients $Y_{n, j_\ell}$ all have different numbers of cusps. Then they are clearly not biholomorphic, but their smooth toroidal compactifications $A_{n, j}^*$ and $A_{n, k}^*$ are biholomorphic by Proposition \ref{prop:IsomOfSTC}. This completes the proof.
\end{pf}
%%%%%%%%%%%%%%%%%%%%

%%%%%%%%%%%%%%%%%%%%
\bibliography{MoreHirzebruch}

\begin{thebibliography}{10}

\bibitem{AMRT}
A.~Ash, D.~Mumford, M.~Rapoport, and Y.~Tai.
\newblock {\em Smooth compactification of locally symmetric varieties}.
\newblock Math. Sci. Press, 1975.
\newblock Lie Groups: History, Frontiers and Applications, Vol. IV.

\bibitem{Goldman}
William~M. Goldman.
\newblock {\em Complex hyperbolic geometry}.
\newblock Oxford Mathematical Monographs. Oxford University Press, 1999.

\bibitem{Hirzebruch}
F.~Hirzebruch.
\newblock Chern numbers of algebraic surfaces: an example.
\newblock {\em Math. Ann.}, 266(3):351--356, 1984.

\bibitem{Holzapfel}
R.-P. Holzapfel.
\newblock Chern numbers of algebraic surfaces---{H}irzebruch's examples are
  {P}icard modular surfaces.
\newblock {\em Math. Nachr.}, 126:255--273, 1986.

\bibitem{Kasparian}
A.~Kasparian.
\newblock Co-abelian toroidal compactifications of torsion free ball quotients.
\newblock arXiv:1201.0099.

\bibitem{Kollar}
J{\'a}nos Koll{\'a}r.
\newblock {\em Shafarevich maps and automorphic forms}.
\newblock M. B. Porter Lectures. Princeton University Press, 1995.

\bibitem{Mok}
Ngaiming Mok.
\newblock Projective algebraicity of minimal compactifications of
  complex-hyperbolic space forms of finite volume.
\newblock In {\em Perspectives in analysis, geometry, and topology}, volume 296
  of {\em Progr. Math.}, pages 331--354. Birkh\"auser, 2012.

\bibitem{Mumford}
D.~Mumford.
\newblock Hirzebruch's proportionality theorem in the noncompact case.
\newblock {\em Invent. Math.}, 42:239--272, 1977.

\bibitem{Tian--Yau}
G.~Tian and S.-T. Yau.
\newblock Existence of {K}\"ahler-{E}instein metrics on complete {K}\"ahler
  manifolds and their applications to algebraic geometry.
\newblock In {\em Mathematical aspects of string theory ({S}an {D}iego,
  {C}alif., 1986)}, volume~1 of {\em Adv. Ser. Math. Phys.}, pages 574--628.
  World Sci. Publishing, 1987.

\bibitem{Wang}
Hsien~Chung Wang.
\newblock Topics on totally discontinuous groups.
\newblock In {\em Symmetric spaces ({S}hort {C}ourses, {W}ashington {U}niv.,
  {S}t. {L}ouis, {M}o., 1969--1970)}, pages 459--487. Pure and Appl. Math.,
  Vol. 8. Dekker, 1972.

\end{thebibliography}

\end{document}